\documentclass[reqno]{gokova}
\usepackage{epsfig,latexsym,graphicx}

\begin{document}
\def\E{\ifmmode{\mathbb E}\else{$\mathbb E$}\fi} 
\def\N{\ifmmode{\mathbb N}\else{$\mathbb N$}\fi} 
\def\R{\ifmmode{\mathbb R}\else{$\mathbb R$}\fi} 
\def\Q{\ifmmode{\mathbb Q}\else{$\mathbb Q$}\fi} 
\def\C{\ifmmode{\mathbb C}\else{$\mathbb C$}\fi} 
\def\H{\ifmmode{\mathbb H}\else{$\mathbb H$}\fi} 
\def\Z{\ifmmode{\mathbb Z}\else{$\mathbb Z$}\fi} 
\def\P{\ifmmode{\mathbb P}\else{$\mathbb P$}\fi} 
\def\T{\ifmmode{\mathbb T}\else{$\mathbb T$}\fi} 
\def\SS{\ifmmode{\mathbb S}\else{$\mathbb S$}\fi} 
\def\DD{\ifmmode{\mathbb D}\else{$\mathbb D$}\fi} 

\renewcommand{\a}{\alpha}
\renewcommand{\b}{\beta}
\renewcommand{\d}{\delta}
\newcommand{\D}{\Delta}
\newcommand{\e}{\varepsilon}
\newcommand{\g}{\gamma}
\newcommand{\G}{\Gamma}
\newcommand{\la}{\lambda}
\newcommand{\La}{\Lambda}
\newcommand{\n}{\nabla}
\newcommand{\var}{\varphi}
\newcommand{\s}{\sigma}
\newcommand{\Sig}{\Sigma}
\renewcommand{\t}{\tau}
\renewcommand{\th}{\theta}
\renewcommand{\O}{\Omega}
\renewcommand{\o}{\omega}
\newcommand{\z}{\zeta}

\newcommand{\ben}{\begin{enumerate}}
\newcommand{\een}{\end{enumerate}}
\newcommand{\be}{\begin{equation}}
\newcommand{\ee}{\end{equation}}
\newcommand{\bea}{\begin{eqnarray}}
\newcommand{\eea}{\end{eqnarray}}
\newcommand{\bc}{\begin{center}}
\newcommand{\ec}{\end{center}}

\newcommand{\IR}{\mbox{I \hspace{-0.2cm}R}}
\newcommand{\IN}{\mbox{I \hspace{-0.2cm}N}}

\newtheorem{thm}{Theorem}[section]
\newtheorem{cor}[thm]{Corollary}
\newtheorem{lem}[thm]{Lemma}
\newtheorem{prop}[thm]{Proposition}
\newtheorem{ax}{Axiom}
\newtheorem{conj}[thm]{Conjecture}

\theoremstyle{definition}
\newtheorem{defn}{Definition}[section]

\theoremstyle{remark}
\newtheorem{rem}{\rm\bfseries{Remark}}[section]
\newtheorem*{notation}{Notation}

\newtheorem{ques}{\rm\bfseries{Question}}[section]
\newtheorem{cons}[rem]{\rm\bfseries{Construction}}
\newtheorem{exm}[rem]{\rm\bfseries{Example}}



\title{Topological quantum field theory and hyperk{\"a}hler geometry}
\author[SAWON]{Justin Sawon
}

\thanks{The author is supported by a Junior Research Fellowship at New
College, Oxford.}

\address{Mathematical Institute \\
         24-29 St Giles' \\
         Oxford OX1 3LB \\
         United Kingdom}
\email{sawon@maths.ox.ac.uk}

\volume{7}

\maketitle

\section{Introduction}

Rozansky and Witten~\cite{rw97} proposed a 3-dimensional
sigma-model whose target space is a hyperk{\"a}hler manifold $X$. For
compact $X$, they conjectured that this theory has an associated
topological quantum field theory (TQFT) with Hilbert spaces given by
certain cohomology groups of $X$. In particular, the Hilbert space
${\mathcal H}_g$ for a genus $g$ Riemann surface should be
$${\mathcal
H}_g:=\bigoplus_q\mathrm{H}^q(X,(\La^{\bullet}T)^{\otimes
g}),$$ 
where we regard $X$ as a complex manifold with respect to some choice
of complex structure compatible with the hyperk{\"a}hler metric
(precisely how these spaces depend on this choice is a subtle
matter). For $X$ a K3 surface, Rozansky and Witten investigated the
cases $g=0$ and $g=1$, and exhibited an action of the mapping class
group in the latter case.

There is a modified TQFT constructed by Murakami and
Ohtsuki~\cite{mo97} using the universal quantum invariant. The Hilbert
spaces in this theory are certain spaces of diagrams, which are
graded modules over a certain commutative ring (we shall make this
precise in due course). This diagrammatic TQFT satisfies a modified
version of the usual TQFT axioms. Let us give some background on the
construction of the Murakami-Ohtsuki TQFT.

The Kontsevich integral~\cite{kontsevich93} was the first construction
of a universal finite-type invariant of links. Le, Murakami, and
Ohtsuki later used the Kontsevich integral to construct an invariant
of closed 3-manifolds, the LMO invariant~\cite{lmo98}. This invariant
is universal among finite-type invariants of integral homology
spheres; this is also the case for rational homology spheres if we use
the Goussarov-Habiro theory of finite-type invariants. The
Murakami-Ohtsuki TQFT is based on a generalization of the LMO
invariant to 3-manifolds with boundary. Hence we believe that, in some
sense, the Murakami-Ohtsuki TQFT should also be regarded as some kind
of universal finite-type object. By applying weight systems to this
``universal finite-type TQFT'' it should be possible to obtain
particular TQFTs.

Rozansky-Witten theory naturally leads to a weight system on graph
cohomology built from a hyperk{\"a}hler manifold. In~\cite{sawon00}
the author constructed a generalization of this weight system to chord
diagrams on circles by adding vector bundles over the hyperk{\"a}hler
manifold (this construction was also discovered independently by
Thompson~\cite{thompson00}). In this article we extend these ideas in
order to apply a ``hyperk{\"a}hler weight system'' to the
Murakami-Ohtsuki TQFT. There are still some difficulties with this
construction (in particular, it is not clear how to apply the weight
system to connected 3-manifolds with disconnected
boundaries). Nevertheless, we are led to a hyperk{\"a}hler TQFT with
the same Hilbert spaces as Rozansky and Witten's. Presumably these are
the same TQFT - given the close relation between the LMO invariant and
the Rozansky-Witten invariant, as investigated by Habegger and
Thompson~\cite{ht99}, this seems like a natural conjecture to make.

Let us outline the contents of this article. We begin in Section 2 by
reviewing the construction of the Murakami-Ohtsuki TQFT. In Section 3
we describe hyperk{\"a}hler manifolds and how they may be used to
construct weight systems. We then apply such a weight system to the
Murakami-Ohtsuki TQFT and describe the hyperk{\"a}hler TQFT so
obtained. In Section 4 we reinterpret the observables of
Rozansky-Witten theory in the context of this TQFT. The main result
here is Proposition~\ref{chern} - in effect we see that all of the
observables can be obtained by pairing vectors from our TQFT with 
cohomology classes. Section 5 is an appendix containing a technical
diagrammatic result required in the construction of the
hyperk{\"a}hler weight system.

Although these ideas should lead to a better understanding of
Rozansky-Witten theory, there are many new questions to explore. For
example, Murakami~\cite{murakami99} described the actions of the
mapping class groups on the Murakami-Ohtsuki TQFT. From this we may
deduce the actions of the mapping class groups on the hyperk{\"a}hler
Hilbert spaces ${\mathcal H}_g$. This will appear in a future
article.

The author would like to express his gratitude to the organizers for
the invitation to speak at the Gokova conference. He would also like
to thank Nigel Hitchin, Thang Le, Christine Lescop, Jun Murakami,
Tomotada Ohtsuki, Justin Roberts, and Simon Willerton for their
comments and many useful discussions.

\section{The Murakami-Ohtsuki TQFT}

In this first section we will describe the modified TQFT of Murakami
and Ohtsuki~\cite{mo97}. To do this we need to understand the
Kontsevich integral of links in $S^3$ and the LMO invariant of
3-manifolds, and we outline the basic ideas behind their
construction. Then we can define the generalizations of these objects
used to construct the TQFT. First we describe the spaces of diagrams
to which these objects belong.

Let $P$ be an oriented 1-dimensional space. In fact, we shall only be
interested in the cases that $P$ is a collection of circles and
trivalent graphs, possible empty. An orientation of a trivalent graph
in this instance is given by an orientation of each edge. A {\em chord
diagram\/} $D$ on $P$ is the union of $P$ and a {\em chord graph\/}
$Q$ - an oriented unitrivalent graph whose univalent vertices lie on
the non-singular part of
$P$. For these unitrivalent graphs, an orientation shall mean an
equivalence class of cyclic orderings of the edges at each trivalent
vertex, with two such being equivalent if they differ at an even
number of vertices. We allow $Q$ to be disconnected, and even to
contain connected components with no univalent vertices. We usually
distinguish $P$ from $Q$ in a chord diagram $D$, but if we do not then
$D$ itself may be regarded as a trivalent graph. Note that those
vertices which occur when a univalent vertex of $Q$ meets $P$ have a
canonical cyclic ordering of edges, induced from the orientation of
$P$: for example, we can take the ordering outgoing edge of $P$,
followed by incoming edge of $P$, followed by the edge belonging to
$Q$. In our diagrams we shall draw $P$ with bold lines to distinguish
it from $Q$\footnote{An alternative convention common in the
literature is to draw $P$ with solid lines and $Q$ with dashed
lines.}.

\begin{figure}[htb]
\includegraphics[scale=0.62]{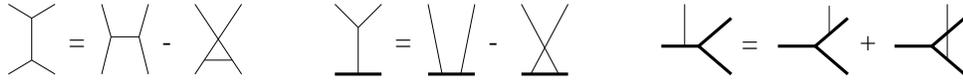}
\caption{The IHX, STU, and branching relations respectively.}
\label{figure1}
\end{figure}

We will consider the space of rational linear combinations of chord
diagrams on $P$. We factor out certain equivalence relations, known as
the AS, IHX, STU, and branching relations. The first of these says
that reversing the orientation of the chord graph $Q$ (i.e.\ reversing
the cyclic ordering of edges at an odd number of trivalent vertices)
is the same as multiplying by $-1$. The remaining relations are as
shown in Figure~\ref{figure1}, with each diagram denoting some part of
a chord diagram. Note that in these and all future diagrams
tetravalent vertices are not vertices at all - they are simply
crossings of edges. Chord diagrams are abstract objects, meaning they
are not embedded in any ambient space. However, since we only have a
sheet of paper on which to draw them we are inevitably hampered by the
dimensional deficiencies of this environment.

When drawn in the plane, we shall assume that our diagrams have the
canonical orientation given by ordering the edges at each vertex
clockwise. In particular, this implies that in the above diagrams the
edges of $P$ (the bold lines) should all be oriented to the right. We
won't always mark the orientation of $P$ on our diagrams as it will
usually be clear from these conventions.

Furthermore, the operation which reverses the orientation of an edge
of $P$ can be made to act on diagrams in a compatible way. Henceforth
we shall generally refrain from discussing questions of orientation of
chord diagrams. The dedicated reader is encouraged to check
compatibility at various stages.

There is a grading on chord diagrams given by half the number of
vertices (univalent and trivalent) of the chord graph $Q$. We denote
the graded completion of the space of chord diagrams on $P$ modulo the
above relations by ${\mathcal A}(P)$. For example, if $P$
is the empty set then ${\mathcal A}(\emptyset)$ is the graded
completion of the space of rational linear combinations of oriented
trivalent graphs modulo the AS and IHX relations, also known as {\em
graph cohomology\/}. This space is a graded commutative ring, with
multiplication given by disjoint union of graphs. In general
${\mathcal A}(P)$ is a graded module over ${\mathcal A}(\emptyset)$,
with the action given by disjoint union of trivalent graphs with
chord diagrams on $P$.

\subsection{The Kontsevich integral}

The Kontsevich integral $Z$ is an invariant of framed oriented links
in $S^3$~\cite{kontsevich93}. For a link $L$ with $l$ components,
$Z(L)$ takes values in the space of chord diagrams on the disjoint
union of $l$ oriented circles
$${\mathcal A}(\coprod_{i=1}^lS^1).$$
The Kontsevich integral is a universal finite-type invariant of links,
which means that all finite-type invariants are obtained by applying
an appropriate weight system to $Z$ (see~\cite{bar-natan95} for
details).

\begin{figure}[htb]
\includegraphics[scale=0.58]{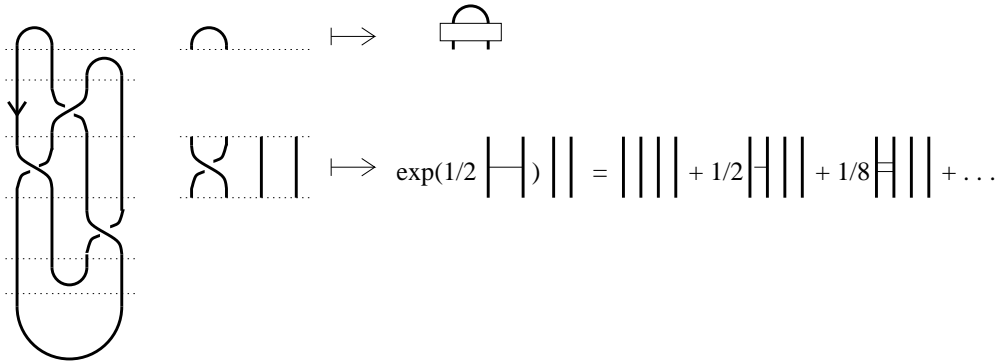}
\caption{Constructing the Kontsevich integral of a trefoil in $S^3$.}
\label{figure2}
\end{figure}

There are various ways to define $Z$; we shall outline the main ideas
behind one way of constructing it. Suppose we have a framed oriented
link $L$ in $S^3$. Then by an ambient isotopy we can ``stretch it
out'' in such a way that taking horizontal slices results in only
fairly simple pieces. Figure~\ref{figure2} shows this for the
trefoil. In fact what we are really using is a projection of $L$ to
the plane, known as a {\em link diagram\/} (not to be confused with a
chord diagram). In such a projection the framing is given by the
blackboard framing. We can define how $Z$ acts on each of these pieces
and then recombine to get complete (linear combinations of) chord
diagrams.

Once we suspect that it may be possible to define $Z$ in this way,
then there is really not much choice in what the individual pieces
should be mapped to, as $Z$ must be invariant under Reidemeister
moves and isotopy of the plane. In particular, $Z$ of the
under-crossing and over-crossing are completely determined. In
Figure~\ref{figure2} we see that the over-crossing is mapped to the
exponential of a single chord multiplied by $1/2$, with additional
strands simply carried through the calculation. For the under-crossing
this factor is $-1/2$. There is some freedom in choosing $Z$ of the
cap and cup, which we have indicated in Figure~\ref{figure2} by
mapping the cap to a box denoting some unspecified collection of
chords. However, these choices only affect the Kontsevich integral up
to an overall normalization. Furthermore, if we assume $Z$ is the
original invariant as defined by Kontsevich~\cite{kontsevich93} then
the normalization is also determined.

\subsection{The LMO invariant}

Le, Murakami, and Ohtsuki~\cite{lmo98} were able to use the Kontsevich
integral to construct an invariant of 3-manifolds, which is a
universal finite-type invariant of rational homology 3-spheres. By a
theorem of Lickorish and Wallace any 3-manifold $M$ can be obtained by
surgery on a framed oriented link $L$ in $S^3$, and we say $L$ {\em
presents\/} $M$. The link $L$ is not uniquely determined, but any
other such link will be obtained from $L$ by a sequence of Kirby
moves. This means that a link invariant which does not distinguish
between links related by Kirby moves can be used to define a
3-manifold invariant.

\begin{figure}[htb]
\includegraphics[scale=0.5]{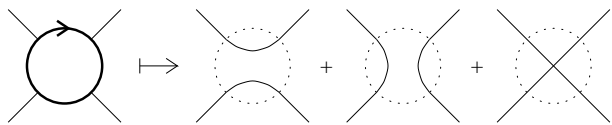}
\caption{Removing a circle using the operation $\iota_2$.}
\label{figure3}
\end{figure}

Given a 3-manifold $M$, Le, Murakami, and Ohtsuki take the Kontsevich
integral $Z(L)$ of some link presenting $M$, and apply to it an
operation $\iota_n$ which removes the circles and replaces them by
part of a unitrivalent graph. If a circle has $m$ legs on it, we
replace it by a an $n$ component tree $T^n_m$ with $m$ legs, and join
these legs to those which were on the circle ($T^n_m$ satisfies a
symmetry property which ensures that the way the legs are joined is
irrelevant). It is not too difficult to show that $T^n_m$ does not 
exist for $m<2n$ and therefore if a chord diagram has fewer than $2n$
legs on some circle it gets mapped to zero by $\iota_n$. When $m=2n$,
$\iota_n$ is given by removing the circle and then summing over all
ways of connecting in pairs the $2n$ legs. Figure~\ref{figure3} shows
this for $n=2$. Le, Murakami, and Ohtsuki further show that any chord
diagram can be expressed as a sum of chord diagrams with at most $2n$
legs on each circle and some additional terms which are killed by
$\iota_n$ (Lemma 3.1 in~\cite{lmo98}). So in effect we only need to
use this `connecting legs in pairs' operation.

After some additional normalization, we get 
$$\O_n(L)=(\mbox{normalizing term})\iota_n(Z(L))\in{\mathcal A}^{\leq n}(\emptyset)$$
where the superscript $\leq n$ is used to denote all diagrams up to
and including degree $n$. The crucial point is that these terms
are invariant under Kirby moves, and hence can be used to build an
invariant of 3-manifolds. In fact, we only need the degree $n$ term in 
$\O_n(L)$ as lower degree terms are contained in
$\O_{n-1}(L)$. Thus we get the LMO invariant
$$\O(M)=1+\sum_{n=1}^{\infty}\O_n(L)^{(n)}\in{\mathcal
A}(\emptyset)$$
of $M$. For rational homology 3-spheres we take a slightly different
normalization $\hat{\O}(M)$, which is given by rescaling terms by
powers of the order of the torsion group ${\mathrm H}_1(M,\Z)$. This
gives us a universal finite-type invariant of rational homology
3-spheres.

As an aside, consider graph cohomology classes of degree one, which
are unique up to scale. They are represented by the {\em theta
graph\/} $\Theta$ for example, and for general $M$ the coefficient of
$\Theta\in{\mathcal A}(\emptyset)$ in $\O(M)$ is Lescop's
generalization of the Casson-Walker invariant~\cite{lescop96}. Thus
the LMO invariant is every bit as powerful as the Casson invariant.

Using these ideas, we can also extend the Kontsevich integral to an
invariant of framed oriented links in arbitrary 3-manifolds as
follows. Suppose $J$ is a link in $M$, and that $L$ presents $M$. Then
there is a link $J^{\prime}$ in $S^3$ such that surgery on $L\subset
S^3$ takes $J^{\prime}$ to $J\subset M$. We take the Kontsevich
integral of $J^{\prime}\cup L$ in $S^3$, and then use the operation
$\iota_n$ as above to remove the circles corresponding to the
components of the link $L$. We can combine the results to get (after
some normalizing) an invariant of $J\subset M$, taking values in
$${\mathcal A}(\coprod_{i=1}^j S^1)$$
where $j$ is the number of components of $J$.

\subsection{The universal quantum invariant for embedded graphs}

\begin{figure}[htb]
\includegraphics[scale=0.7]{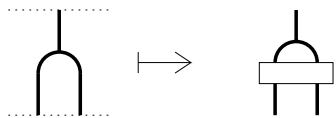}
\caption{The Kontsevich integral of a trivalent vertex.}
\label{figure4}
\end{figure}

Now suppose that instead of a link in $S^3$, we have an embedded
framed oriented trivalent graph. Call the abstract graph $\G$, and
denote its embedding in $S^3$ by $G$. Then one can construct an
invariant of $G\subset S^3$ taking values in ${\mathcal A}(\G)$ in the
same way as we constructed the Kontsevich integral above. The only new
feature is that we need to define the invariant on a piece of $G$
containing a trivalent vertex, as shown in Figure~\ref{figure4}. 

\begin{figure}[htb]
\includegraphics[scale=0.7]{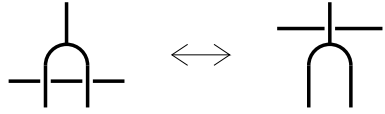}
\caption{Reidemeister move for trivalent graphs.}
\label{figure5}
\end{figure}

As before the box denotes some collection of chords, uniquely
determined by requiring that the overall invariant thus obtained is
invariant under isotopy of the plane and by the Reidemeister-type move
shown in Figure~\ref{figure5}.

This is precisely how Murakami and Ohtsuki~\cite{mo97} define what
they call the {\em universal Vassiliev-Kontsevich invariant\/} of the
graph $G$ in $S^3$, and denote by $\hat{Z}(G)$. As with the Kontsevich 
integral, we can extend this to an invariant of framed oriented
trivalent graphs embedded in an arbitrary 3-manifold by combining with
the LMO invariant. Murakami and Ohtsuki call this generalization the
{\em universal quantum invariant\/} of the graph $G$ in $M$, denoted
$\O(M,G)$. When $M$ is a rational homology 3-sphere we once again have
a slightly different normalization, denoted $\hat{\O}(M,G)$.

\subsection{The LMO invariant for 3-manifolds with boundary}

So far we have only discussed closed 3-manifolds. Now let $M$ be a
3-manifold with boundary $\partial M$. The following construction
works equally well when $\partial M$ has several connected components,
but for ease of exposition we shall assume that $\partial M$ consists
of a single Riemann surface of genus $g$.

\begin{figure}[htb]
\includegraphics[scale=0.6]{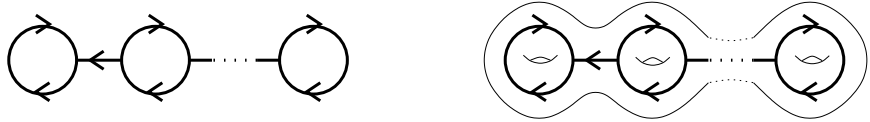}
\caption{The chain graph $\G_g$ and embedded in a neighbourhood $N(G)$.}
\label{figure6}
\end{figure}

Let $\G_g$ be the {\em chain graph\/} with $g$ loops as shown in
Figure~\ref{figure6}. It is given the blackboard framing and its
orientation is as shown. If $\G_g$ is embedded in some ambient
3-dimensional space then its neighbourhood $N(G)$ is a 3-manifold
whose boundary is a genus $g$ surface $\Sig_g$.

Suppose we are given an identification of $\Sig_g$ with $\partial
M$. This requires choosing meridians and longitudes on $\partial M$
and identifying them with standard ones on $\Sig_g$. Once we've done
this, we can glue the neighbourhood $N(G)$ of the embedded $\G_g$ into
$M$, identifying the boundaries $\Sig_g$ and $\partial M$. The result
is a closed 3-manifold $\hat{M}$, which contains $\G_g$ as an embedded
graph $G$. We can regard the universal quantum invariant
$\O(\hat{M},G)$ (or $\hat{\O}(\hat{M},G)$ when $\hat{M}$ is a rational
homology 3-sphere) of $G\subset\hat{M}$ as the LMO invariant of the
3-manifold $M$ with boundary. In the next subsection we shall use this
to construct a modified TQFT.

\subsection{The modified TQFT axioms}

A 3-dimensional TQFT is a functor from the category of oriented
3-cobordisms to a category of modules over a commutative ring $k$. In
concrete terms, we associate to each Riemann surface $\Sig$ a module 
$V(\Sig)$ and to each 3-manifold $M$ with boundary $\Sig=\partial M$
an element $Z(M)\in V(\Sig)$. We require the following axioms to
hold.
\begin{itemize}
\item[(A1)] For surfaces $\Sig_1$ and $\Sig_2$,
$V(\Sig_1\sqcup\Sig_2)=V(\Sig_1)\otimes V(\Sig_2)$\footnote{The
subscripts here are labels and do not indicate surfaces of genus one
and two. Despite this clash of notation, it should be clear what is
meant from the context.}.
\item[(A2)] Reversing the orientation of $\Sig$ gives the dual module
$V(-\Sig)=V(\Sig)^*$.
\item[(A3)] For the empty surface $V(\emptyset)=k$. In particular,
$Z(M)\in k$ for a closed 3-manifold.
\item[(A4)] Suppose $M_1$ and $M_2$ have boundaries $\partial
M_1=\Sig\sqcup\Sig_1$ and $\partial M_2=(-\Sig)\sqcup\Sig_2$
respectively. We can glue $M_1$ and $M_2$ along $\Sig$ to get
$M=M_1\cup_{\Sig}M_2$. Then $Z(M)=\langle
Z(M_1),Z(M_2)\rangle_{\Sig}$, where 
$$\langle\phantom{M},\phantom{M}\rangle_{\Sig}:V(\Sig)\otimes
V(\Sig_1)\otimes V(\Sig)^*\otimes V(\Sig_2)\longrightarrow V(\Sig_1)\otimes
V(\Sig_2)$$
denotes the contraction mapping.
\end{itemize}

In~\cite{mo97} Murakami and Ohtsuki construct the following modified
TQFT. Firstly, we take a sub-category of the category of 3-cobordisms 
such that all the closed 3-manifolds involved will be rational
homology 3-spheres (this amounts to some condition on the homology
groups of the cobordisms). For a genus $g$ Riemann surface $\Sig_g$ we
define
$$V(\Sig_g):={\mathcal A}(\G_g),$$
which is a module over the commutative ring ${\mathcal
A}(\emptyset)$. For disconnected Riemann surfaces we define
$$V(\Sig_1\sqcup\Sig_2):={\mathcal A}(\G_1\sqcup\G_2),$$
and in this way $V(\Sig)$ is defined for all Riemann surfaces $\Sig$.

Let $M$ be a 3-manifold with boundary $\partial M$. Suppose the
boundary is isomorphic to the genus $g$ Riemann surface $\Sig_g$. Then
$V(\partial M)\cong V(\Sig_g)$, and more specifically, {\em for each
choice of a set of longitudes and meridians on\/} $\partial M$ there
is an isomorphism between $V(\partial M)$ and $V(\Sig_g)={\mathcal
A}(\G_g)$ (we assume $\Sig_g$ comes equipped with a standard set of
longitudes and meridians). As in the last subsection, we can glue a
neighbourhood of a graph to $M$ and obtain a closed 3-manifold
$\hat{M}$ containing $\G_g$ as an embedded graph $G$. Because of our
choice of sub-category of 3-cobordisms, $\hat{M}$ is a rational
homology 3-sphere, and hence we can define
$$Z(M):=\hat{\O}(\hat{M},G)\in {\mathcal A}(\G_g)\cong V(\partial M).$$
It is important to remember that the isomorphism on the right depends
on a choice of meridians and longitudes for $\partial M$. This
definition clearly extends to arbitrary (disconnected) Riemann
surfaces $\partial M$.

The first thing to observe is that axiom (A1) is not satisfied;
instead we have
\begin{itemize}
\item[(A1)$^{\prime}$] There is an inclusion $V(\Sig_1)\otimes
V(\Sig_2)\hookrightarrow V(\Sig_1\sqcup\Sig_2)$. 
\end{itemize}
This inclusion is given by disjoint union of chord diagrams on $\G_1$
and $\G_2$. General chord diagrams on $\G_1\sqcup\G_2$ may contain
chord graphs which connect $\G_1$ and $\G_2$, whereas the image
of the above inclusion only contains chord diagrams for which $\G_1$
and $\G_2$ are in distinct connected components.

\begin{figure}[htb]
\includegraphics[scale=0.6]{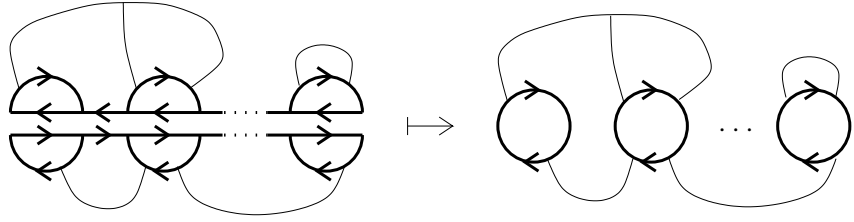}
\caption{Pairing of chord diagrams on $\G_g$ and $-\G_g$.}
\label{figure7}
\end{figure}

Axiom (A2) is replaced by
\begin{itemize}
\item[(A2)$^{\prime}$] There is a pairing 
$$\langle\phantom{M},\phantom{M}\rangle_{\Sig}:V(\Sig)\otimes
V(-\Sig)\longrightarrow {\mathcal A}(\emptyset).$$
\end{itemize}
We now describe this pairing; it suffices to do so for a genus $g$
surface $\Sig_g$. In other words, we have a pair of chord diagrams on
$\G_g$ and $-\G_g$ respectively, which we align as shown in
Figure~\ref{figure7}. Note that the orientations may be arranged as
shown by applying the operation which reverses the orientations of
edges of $\G$ and acts on chord diagrams in a compatible way.

Using the branching relations, we may assume that our chord diagrams
have no chords attached to the straight edges of the graphs, i.e.\ all
legs of the unitrivalent graphs end on the curved edges of $\G_g$ and
$-\G_g$. Now we remove the straight edges and connect the curved
edges to form $g$ oriented circles. Finally, we remove the circles
using the same argument as for the LMO invariant. Thus we are left
with an element of ${\mathcal A}(\emptyset)$, which is precisely what
we wanted. Note that by identifying $V(\Sig)$ with the space of chord
diagrams on $\G$ we have implicitly assumed that we have chosen a set
of meridians and longitudes on $\Sig$. The pairing depends on this
choice.

Axiom (A3) needs no modification. The commutative ring $k$ is
${\mathcal A}(\emptyset)$. Multiplication in ${\mathcal A}(\emptyset)$
and the action on ${\mathcal A}(\G)$ are given by disjoint union of
diagrams.

Suppose we have 3-manifolds $M_1$ and $M_2$ with boundaries
$\partial M_1\cong\Sig$ and $\partial M_2\cong -\Sig$
respectively. Then we can pair $Z(M_1)\in {\mathcal A}(\G)$ with
$Z(M_2)\in {\mathcal A}(-\G)$ as described above, to get an element
$$\langle Z(M_1),Z(M_2)\rangle_{\Sig}\in{\mathcal A}(\emptyset).$$
As usual, this requires a choice of longitudes and meridians on
$\partial M_1$ and $\partial M_2$, so that we can identify them with
$\Sig$ and $-\Sig$ respectively. This induces a homeomorphism $f$ of
the boundaries of our 3-manifolds, with which we can glue $M_1$ and
$M_2$ to get a closed 3-manifold $M=M_1\cup_fM_2$. Murakami and
Ohtsuki~\cite{mo97} prove that $Z(M)\in{\mathcal A}(\emptyset)$ is
essentially the same as the pairing of $Z(M_1)$ and $Z(M_2)$.

More generally, $M_1$ and $M_2$ may have additional boundary
components $\Sig_1$ and $\Sig_2$, respectively. Then we have the
following modification of axiom (A4).
\begin{itemize}
\item[(A4)$^{\prime}$] Under the above hypotheses on $M_1$ and $M_2$, we
have
$$Z(M)=(\mbox{normalizing terms})(\mbox{framing anomaly
correction})\langle Z(M_1),Z(M_2)\rangle_{\Sig}$$
as elements of $V(\Sig_1\sqcup\Sig_2)$.
\end{itemize}
As usual the normalizing terms involve the orders of the integral
homology groups of $\hat{M_1}$, $\hat{M_2}$, and $\hat{M}$ (which are
finite groups since these are rational homology 3-spheres). The {\em
framing anomaly\/} is a subject we'd rather avoid here. Murakami and
Ohtsuki expect that this correction term can be removed by defining a
suitable framing of 3-manifolds (the interested reader should
consult~\cite{mo97}).

The key to proving that the Murakami-Ohtsuki TQFT satisfies axiom
(A4)$^{\prime}$ is the following result (Lemma 4.3 in~\cite{mo97}).
\begin{lem}
Let $M_1$ and $M_2$ be as given above. Suppose that $\hat{M}_1$ and
$\hat{M}_2$ are given by surgery on links $L_1$ and $L_2$ in $S^3$,
respectively. Recall that these closed 3-manifolds contain embedded
graphs $G_1$ and $G_2$. We use the same notation to denote the graphs 
in $S^3$ which become $G_1$ and $G_2$ after the surgeries on $L_1$ and
$L_2$. Let $L$ be the link in $S^3$ obtained from $L_1\cup G_1$ and
$L_2\cup G_2$ by removing the straight edges of $G_1$ and $G_2$, and
connecting the curved edges to form embedded circles (we may need to
first isotope the trivalent graphs into adjacent positions). Then $M$
is homeomorphic to the 3-manifold given by surgery on $L$.
\end{lem}
\noindent
In calculating the LMO invariant $Z(M)$ of $M$ we perform the circle
removing operation on $L$, so we have really done the same operations
as we would to calculate the pairing $\langle
Z(M_1),Z(M_2)\rangle_{\Sig}$. Thus axiom (A4)$^{\prime}$ is
satisfied. 

Note that for a closed 3-manifold $M$, $Z(M)$ is the LMO invariant
$\hat{\O}(M)\in{\mathcal A}(\emptyset)$, with the rational homology
3-sphere normalization. If $M$ has boundary $\partial M\cong S^2$, we
also take $Z(M)\in{\mathcal A}(\emptyset)$. In other words, the genus
zero chain graph $\G_0$ is the empty graph, but we think of a
neighbourhood of $\G_0$ as being a solid ball with boundary $S^2$. The
reader can think of this as simply a convention, though in some ways
it is the only sensible choice.

\section{Hyperk{\"a}hler geometry}

It is common in the theory of knot and 3-manifold invariants to apply
Lie algebra weight systems to the universal diagrammatic invariants -
this is one way of generating quantum invariants. Indeed the words
`universal quantum invariant' mean that all quantum invariants arise
in this way. For example, Jones' polynomial arises by applying an
$\mathfrak{su}(2)$ weight system to the Kontsevich integral. In this
section we will describe how a {\em hyperk{\"a}hler weight system\/}
can be applied to the Murakami-Ohtsuki TQFT. In some ways this is a
continuation of the author's work in~\cite{sawon00}, where the genus
one case was described (see also Thompson~\cite{thompson00}). The
argument presented here works only for 3-manifolds with connected
boundaries, and in this sense our programme is still
incomplete. Nevertheless, we expect that it should be possible to
extend the results to the disconnected case.

\subsection{Hyperk{\"a}hler and holomorphic symplectic manifolds}

Let $X$ be a compact\footnote{To construct Rozansky-Witten invariants, 
as in~\cite{rw97}, the assumption that the manifold is compact may be
dropped provided the appropriate asymptotic decay conditions are
satisfied by the curvature. For the construction of the TQFT we assume
$X$ is compact, and are uncertain whether this can be generalized to
non-compact manifolds.} irreducible\footnote{The theory easily extends
to reducible manifolds.} hyperk{\"a}hler manifold of
real-dimension $4k$. In other words, there is a metric on $X$ whose
Levi-Civita connection has holonomy $\mathrm{Sp}(k)$. Such a manifold
admits a triple of complex structures $I$, $J$, and $K$, which act
like the quaternions on the tangent bundle $TX$. A {\em hyperk{\"a}hler
metric\/} is a metric $g$ which is K{\"a}hlerian with respect to all
of these complex structures. Let us call the corresponding K{\"a}hler
forms $\o_1$, $\o_2$, and $\o_3$ respectively. 

There is no canonical choice of complex structure on $X$ compatible
with the metric, but since we wish to use the methods of complex
geometry we shall fix a structure $I$ and henceforth regard $X$ as a
complex manifold of complex-dimension $2k$. Then $\o_2$ and $\o_3$ can
be combined to give a two-form
$$\o:=\o_2+i\o_3\in\mathrm{H}^0(X,\La^2T^*)$$
which is holomorphic with respect to $I$ ($T$ now denotes the
holomorphic tangent bundle). 

Rozansky-Witten theory~\cite{rw97} was originally built around
hyperk{\"a}hler manifolds, though Kontsevich~\cite{kontsevich99i} and
Kapranov~\cite{kapranov99} later showed that all one requires is a
complex manifold with a holomorphic symplectic form, otherwise known
as a {\em holomorphic symplectic manifold\/}. In particular, if it is
K{\"a}hler then it must be a hyperk{\"a}hler manifold (as we are
working in the compact setting), but there are also non-K{\"a}hler
examples due to Guan~\cite{guan95}. However, we suspect that there may
be further interesting properties for hyperk{\"a}hler $X$ which are
not present for general holomorphic symplectic manifolds. For example,
it may be interesting to observe the way that the Hilbert spaces of
the TQFT change under a variation of the compatible complex structure
on $X$. So although in this article the reader may assume that $X$ is
merely a holomorphic symplectic manifold, we will stick to the
original ``hyperk{\"a}hler'' terminology nonetheless.

The key to Kapranov's version of the Rozansky-Witten weight system is
to use the Atiyah class~\cite{atiyah57} instead of the curvature of
the hyperk{\"a}hler manifold. The Atiyah class 
$$\a_E\in\mathrm{H}^1(X,T^*\otimes\mathrm{End}E)$$
of a complex vector bundle $E$ on $X$ is the obstruction to the
existence of a global holomorphic connection on $X$. The case $E=T$
will be the only one of interest to us, for which the Atiyah class
lies in
$$\mathrm{H}^1(X,T^*\otimes T^*\otimes T).$$
Observe that we can identify $T$ and $T^*$ using the holomorphic
symplectic form $\o$. Kapranov showed that the element $\a_T$ is
totally symmetric (Proposition 5.1.1 in~\cite{kapranov99}), i.e.\ it
lies in
$$\mathrm{H}^1(X,\mathrm{Sym}^3T^*).$$
In local complex coordinates, we write $\a_{ijk}$ for $\a_T$, where
the subscripts refer to the three copies of $T^*$.

The last thing we shall need in the construction is the dual of the
holomorphic symplectic form
$$\tilde{\o}\in\mathrm{H}^0(X,\La^2T).$$
In local complex coordinates, $\tilde{\o}$ has matrix $\o^{ij}$. Note
that this is {\em minus\/} the inverse of the matrix $\o_{ij}$ of
$\o$.

\subsection{Some diagrammatic preliminaries}

Suppose we have a 3-manifold $M$ with boundary $\partial M$ a genus
$g$ Riemann surface. Then in the Murakami-Ohtsuki TQFT we get, after a
choice of longitudes and meridians in $\partial M$, an element
$$Z(M)\in\mathcal{A}(\G_g).$$
By applying a hyperk{\"a}hler weight system to $Z(M)$ we would like to
obtain something in
$${\mathcal H}_g=\bigoplus_q\mathrm{H}^q(X,(\La^{\bullet}T)^{\otimes
g}).$$ 
First we need to rewrite elements of $\mathcal{A}(\G_g)$ in a nicer
way.

A {\em marked unitrivalent graph\/} is an oriented unitrivalent graph
$D$ whose univalent vertices (or {\em legs\/}) are labelled by the
integers $1$ to $g$. An orientation of such a graph is an equivalence
class of cyclic orderings of the edges at each trivalent vertex, with
two such being equivalent if they differ at an even number of
vertices. Note that more than one leg may be labelled by the same
integer, and we do not need to use all labels. The graphs may also be
disconnected, and may contain connected components with no legs.

We will consider the space of rational linear combinations of marked
unitrivalent graphs modulo the AS and IHX relations. This space is
graded by half the number of vertices (univalent and trivalent) of a
unitrivalent graph $D$, and we denote the graded completion by
$\mathcal{B}_g$.

Let $I$ be the interval, and let
$$\chi :\mathcal{B}_g\longrightarrow\mathcal{A}(\coprod_{i=1}^gI)$$
be the map given by averaging over all ways of joining the legs of a
marked unitrivalent graph to the intervals such that the legs labelled
by $j\in\{1,\ldots,g\}$ are joined to the $j^{\mathrm{th}}$ interval.

\begin{prop}
The map $\chi$ is an isomorphism of $\mathcal{A}(\emptyset)$-modules.
\end{prop}
\noindent
{\bf Proof:} This result is well-known among knot-theorists, and is
described for $g=1$ in Bar-Natan~\cite{bar-natan95}\footnote{Note that
Bar-Natan uses the map given by summing over all ways of joining the
legs to the interval, whereas we use the average.} and for general $g$
in Bar-Natan et al~\cite{bgrtI99} (Definition 2.7). In the former
case, the argument relies on using the IHX relations to rewrite an
element of $\mathcal{B}:=\mathcal{B}_1$ in such a way that its legs
are arranged `symmetrically'. The map $\chi$ is then simply given by
gluing the legs to the interval and is an isomorphism. Since this is a
`local' operation, it generalizes to arbitrary $g$. In other words, we
can rewrite an element of $\mathcal{B}_g$ in such a way that its legs
are arranged in $g$ symmetric collections, and the proposition
follows.\hfill$\Box$

\begin{figure}[htb]
\includegraphics[scale=0.6]{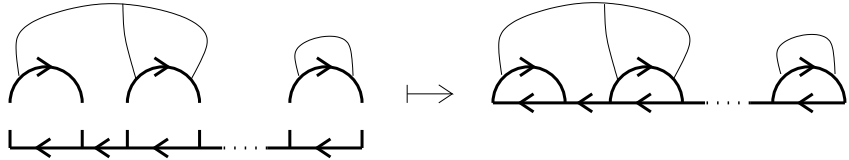}
\caption{Attaching a tree to a collection of intervals.}
\label{figure8}
\end{figure}

\vspace*{5mm}
There is a surjective map
$$\mathcal{A}(\coprod_{i=1}^gI)\longrightarrow\mathcal{A}(\coprod_{i=1}^gS^1)$$
given by closing up the intervals into circles. This is only an
isomorphism when $g=1$ (Theorem 8 in
Bar-Natan~\cite{bar-natan95}; see also Lemma 3.17 in
Thurston~\cite{thurston00} and the discussion there). However, this is
not the map we wish to use. Instead we wish to attach a tree to the
intervals as shown in Figure~\ref{figure8}, resulting in a chain graph
$\G_g$. This gives us a map
$$\rho :\mathcal{A}(\coprod_{i=1}^gI)\longrightarrow\mathcal{A}(\G_g).$$

\begin{prop}
The map $\rho$ is an isomorphism of $\mathcal{A}(\emptyset)$-modules.
\label{rho}
\end{prop}
\noindent
{\bf Proof:} We postpone the proof to an appendix.\hfill$\Box$

\vspace*{5mm}
Composing $\chi$ and $\rho$ gives us an isomorphism
$$\tau:=\rho\circ\chi :\mathcal{B}_g\longrightarrow\mathcal{A}(\G_g).$$
This means that given a chord diagram $D$ on the chain graph $\G_g$,
we can rewrite it in an equivalent way as a marked unitrivalent graph,
namely $\tau^{-1}(D)$.

\subsection{The hyperk{\"a}hler weight system}

Recall that $\mathcal{B}_g$ is graded by half the number of
vertices. It also has a multi-grading given by the number of legs with
label 1, with label 2, etc. A unitrivalent graph with no legs is
simply a trivalent graph, and thus
$\mathcal{B}_0=\mathcal{A}(\emptyset)$. Our aim in this subsection is
to show the following.

\begin{prop}
There is a homomorphism of vector spaces
$$W_X:\mathcal{B}_g\longrightarrow\mathcal{H}_g:=\bigoplus_q\mathrm{H}^q(X,(\La^{\bullet}T)^{\otimes g}).$$
If a marked unitrivalent graph $D$ has $q$ trivalent vertices and
$l_i$ legs labelled by $i\in\{1,\ldots,g\}$, then $W_X(D)$ lies in
$$\mathrm{H}^q(X,\La^{l_1}T\otimes\cdots\otimes\La^{l_g}T).$$
The map 
$$W_X|_{\mathcal{B}_0=\mathcal{A}(\emptyset)}:\mathcal{A}(\emptyset)\longrightarrow\mathcal{H}:=\mathcal{H}_0=\bigoplus_q\mathrm{H}^q(X,\mathcal{O}_X)$$
is a homomorphism of commutative rings, where the product in
$\mathcal{H}$ is cup-product. In general the map $W_X$ intertwines
the $\mathcal{A}(\emptyset)$-module and $\mathcal{H}$-module
structures on $\mathcal{B}_g$ and $\mathcal{H}_g$ respectively.
\end{prop}
\noindent
{\bf Proof:} We will give a definition of $W_X$, after which the rest
of the statements in the proposition will follow. We begin with
$g=0$. This is the situation of the original Rozansky-Witten
invariants~\cite{rw97}. We review instead the approach of
Kapranov~\cite{kapranov99} which leads immediately to cohomology
classes.

Suppose $D$ is a trivalent graph with $q$ vertices (in this
case $q$ must be even). We place a copy of the Atiyah class $\alpha_T$
at each vertex of $D$ and a copy of the dual of the holomorphic
symplectic form $\tilde{\o}$ at each edge. If $\alpha_{ijk}$ has been
placed at some vertex, we label the outgoing edges by $i$, $j$, and
$k$. Similarly, if $\o^{ij}$ has been placed on some edge, we label
the ends of the edge by $i$ and $j$\footnote{These labellings must be
done in a way compatible with the orientation of $D$. This is
discussed at length in Kapranov~\cite{kapranov99}, and reproduced in
Hitchin and Sawon~\cite{hs99}. We do not wish to repeat the argument
for a third time, so the reader is urged to consult those articles
regarding questions of orientation.}. We then contract all indices and
take the cup product of the $\alpha_T$s. This results in an element
$$W_X(D)\in\mathrm{H}^q(X,\mathcal{O}_X)$$
which is precisely what we wanted. That this construction is
compatible with the AS and IHX relations follows (respectively) from
careful consideration of the orientation of $D$ and from an identity
satisfied by the Atiyah class. See Kapranov~\cite{kapranov99} or
Hitchin and Sawon~\cite{hs99} for details.

The general case is much the same. Let $D$ be a marked
unitrivalent graph with $q$ trivalent vertices and $l_i$ legs labelled
by $i\in\{1,\ldots,g\}$. Since we have the isomorphism
$$\chi :\mathcal{B}_g\longrightarrow\mathcal{A}(\coprod_{i=1}^gI),$$
we may assume that the legs of $D$ are already arranged into $g$
symmetric collections. As before, we place a copy of $\alpha_T$ at
each trivalent vertex and a copy of $\tilde{\o}$ at each edge, and
label outgoing edges and ends of edges by indices. Contracting indices
and taking the cup product of the $\alpha_T$s gives us an element of
$$\mathrm{H}^q(X,T^{\otimes (l_1+\ldots +l_g)}),$$
where the copies of $T$ occur because there are uncontracted indices
labelling the legs. The fact that the legs are arranged into $g$
symmetric collections means that we actually get an element
$$W_X(D)\in\mathrm{H}^q(X,\La^{l_1}T\otimes\cdots\otimes\La^{l_g}T)\footnote{The
fact that the legs are arranged {\em symmetrically\/} but give us
a cohomology class with values in an {\em anti-symmetric\/} bundle,
namely $\La^{l_1}T\otimes\cdots\otimes\La^{l_g}T$, is a manifestation
of the general `reversal of statistics' from Bose-Einstein to
Fermi-Dirac which is inherent in Rozansky-Witten theory. The same
phenomena occurs with the wheels in Hitchin and Sawon~\cite{hs99}, and
is the reason why Kapranov considers desuspensions of operads
in~\cite{kapranov99}. If you think this is some kind of skullduggery,
then check it directly by considering the orientations!}.$$
Compatibility with the AS and IHX relations follow as
before.\hfill$\Box$

\subsection{The hyperk{\"a}hler TQFT}

We wish to define a TQFT by applying the hyperk{\"a}hler weight system
of the previous subsection to the Murakami-Ohtsuki TQFT. Recall that
the Hilbert spaces of the latter are modules over the commutative ring
$\mathcal{A}(\emptyset)$. Elements of this ring are mapped to the
space
$$\mathcal{H}=\bigoplus_q \mathrm{H}^q(X,\mathcal{O}_X)$$
which will be the commutative ring of our hyperk{\"a}hler TQFT. Note
that for irreducible hyperk{\"a}hler manifolds the cohomology groups
$\mathrm{H}^q(X,\mathcal{O}_X)$ vanish for odd $q$ and are one
dimensional and generated by
$$[\bar{\o}^l]\in\mathrm{H}^{0,2l}_{\bar{\partial}}(X)=\mathrm{H}^{2l}(X,\mathcal{O}_X)$$
for even $q=2l$. Up to scale, there is a unique graph cohomology class
of degree one, represented by the graph $\Theta$ for example. When $X$
is irreducible
$$W_X(\Theta)=\beta_{\Theta}[\bar{\o}]\in\mathrm{H}^2(X,\mathcal{O}_X)$$
for some scalar $\beta_{\Theta}$. This scalar is proportional to the
$\mathcal{L}^2$-norm of the curvature of $X$ and is therefore non-zero
(see Equation 9 in Hitchin and Sawon~\cite{hs99}). Similarly, for the
disjoint union of $l$ copies of $\Theta$ we have
$$W_X(\Theta^l)=\beta_{\Theta}^l[\bar{\o}^l]\in\mathrm{H}^{2l}(X,\mathcal{O}_X).$$
Thus $W_X(\Theta)$ generates the commutative ring $\mathcal{H}$.

In the original paper of Rozansky and Witten~\cite{rw97} there was a 
further map from $\mathcal{H}$ to the real numbers. This was given by
taking trivalent graphs of degree $k$ where the real-dimension of $X$
is $4k$; $W_X$ applied to such a graph gives us an element of
$\mathrm{H}^{2k}(X,\mathcal{O}_X)$, and taking the Serre duality
pairing with the generator $[\o^{2k}]$ of
$\mathrm{H}^0(X,\La^{2k}T^*)$ gives us a number (which is real after
the appropriate normalization). If we begin with a trivalent graph of
degree less than $k$, it is still possible to get a real number by
including some {\em observables\/}. Basically this involves
multiplying by some power of $[\bar{\o}]$, or taking the disjoint
union with some copies of $\Theta$ (these are equivalent for $X$ compact
and irreducible). We shall say more about these and other observables
in the next section. From the TQFT perspective, however, it makes more
sense to work over the commutative ring $\mathcal{H}$ in order to
preserve more of the original structure.

Recall that the Hilbert spaces of the Murakami-Ohtsuki TQFT are the
spaces $\mathcal{A}(\G_g)$. Applying the isomorphism $\tau^{-1}$ takes
us to $\mathcal{B}_g$, to which we apply the hyperk{\"a}hler weight
system $W_X$, taking us into
$$\mathcal{H}_g=\bigoplus_q\mathrm{H}^q(X,(\La^{\bullet}T)^{\otimes g}).$$
We must be careful here. Whereas we can ensure that the map
$$W_X:\mathcal{A}(\emptyset)\longrightarrow\mathcal{H}$$
is onto, this is never likely to be the case for
$$W_X\circ\tau^{-1}:\mathcal{A}(\G_g)\longrightarrow\mathcal{H}_g.$$
This suggests that the Hilbert spaces of our hyperk{\"a}hler TQFT
should be subspaces of $\mathcal{H}_g$, though exactly how to define
these subspaces is still unclear.

Putting aside this problem, let us just assume the Hilbert spaces are
the $\mathcal{H}_g$s themselves. Then a 3-manifold $M$ with boundary
$\partial M\cong\Sig_g$ gives rise to an element
$$Z(M)\in\mathcal{A}(\G_g)$$
after a choice of longitudes and meridians in $\partial M$. Applying
$W_X\circ\tau^{-1}$ takes us into $\mathcal{H}_g$ and completes the
definition of our TQFT.

Note that we assumed $M$ has a single connected boundary
component. The application of the weight system $W_X$ in the more
general case remains problematical. Presumably we would like axiom
(A1)
$$V(\Sig_1\sqcup\Sig_2)=V(\Sig_1)\otimes V(\Sig_2)$$
to be satisfied. For example
$$V(\Sig_{g_1}\sqcup\Sig_{g_2})=\mathcal{H}_{g_1}\otimes_{\mathcal{H}}\mathcal{H}_{g_2}$$
and every Hilbert space will be a tensor product of $\mathcal{H}_g$s,
tensored over $\mathcal{H}$ since these are
$\mathcal{H}$-modules. However, the Murakami-Ohtsuki TQFT satisfies 
the modified axiom (A1)$^{\prime}$
$$V(\Sig_1)\otimes V(\Sig_2)\hookrightarrow V(\Sig_1\sqcup\Sig_2)$$
i.e.\
$$\mathcal{A}(\G_1)\otimes_{\mathcal{A}(\emptyset)}\mathcal{A}(\G_2)\hookrightarrow\mathcal{A}(\G_1\sqcup\G_2).$$
Although the left hand side sits inside the right as a direct summand,
projecting onto it results in a significant loss of information. On
the other hand, there does not appear to be an isomorphism from
$\mathcal{A}(\G_1\sqcup\G_2)$ to the space of marked unitrivalent
graphs (or some similar space) which generalizes $\tau^{-1}$. Thus it
is not clear how to extend the hyperk{\"a}hler weight system $W_X$ in
this case.

Next let us consider the pairing of axiom (A2)$^{\prime}$. Let $D_1$
and $D_2$ be chord diagrams on $\G_g$ and $-\G_g$
respectively. Note that a chord diagram on $-\G_g$ is the same as a
chord diagram on $\G_g$, up to a sign, so we may as well take
$D_2\in\mathcal{A}(\G_g)$. We described the pairing
$$\langle D_1,D_2\rangle_{\Sig_g}\in\mathcal{A}(\emptyset)$$
in Subsection 2.5. Applying the hyperk{\"a}hler weight system directly
$$\langle W_X\circ\tau^{-1}(D_1),W_X\circ\tau^{-1}(D_2)\rangle_{\mathcal{H}_g}=W_X(\langle D_1,D_2\rangle_{\Sig_g})\in\mathcal{H}$$
gives us an implicit description of the pairing in the hyperk{\"a}hler
TQFT. This is somewhat unsatisfactory. We'd like to have an explicit
description of the pairings in the spaces $\mathcal{H}_g$. In the
hyperk{\"a}hler context, the operation of summing over all pairings of
legs becomes a combination of wedge product of $\La^{\bullet}T$s and
fiberwise convolution, i.e.\
$$\La^{2m}T^*\otimes\La^{2m}T\longrightarrow\mathcal{O}_X,$$
with a power of the holomorphic symplectic form $\omega$, which gives
a section of $\La^{2m}T^*$ (an example of such a calculation was
carried out in detail in Hitchin and Sawon~\cite{hs99}). Combining
this with cup-product of cohomology classes no doubt gives the pairing
on $\mathcal{H}_g$, though unfortunately this is rather difficult to
see from the implicit description given above. The problem is that in
order to reduce the LMO `circle removing operation' (see Subsection
2.2) to a sum over all pairings of legs, we first need to rewrite the 
chord diagrams $D_1$ and $D_2$ in some canonical way. On the other
hand, we need to apply the isomorphism $\tau^{-1}$ before acting with
the hyperk{\"a}hler weight system $W_X$. This is also equivalent to
choosing some canonical way of writing $D_1$ and $D_2$, and these two
canonical representations are unlikely to be compatible.

The third axiom (A3) simply tells us that closed 3-manifolds will
gives us $\mathcal{H}$-valued invariants. These are the original
Rozansky-Witten invariants of~\cite{rw97}, after including observables
where necessary and mapping $\mathcal{H}$ to the real numbers, as
described earlier in this subsection. The final axiom (A4)$^{\prime}$
should still be satisfied, though as we have yet to fully extend our 
hyperk{\"a}hler TQFT to all 3-manifolds we will postpone our
investigations to a future article.

\section{Observables}

There are three kinds of observables occurring in the Rozansky-Witten
theory. The first was introduced by Rozansky and
Witten~\cite{rw97} and the second and third by
Thompson in~\cite{thompson00} and~\cite{thompson99} respectively. In
this subsection we review these and try to place them within the
context of our hyperk{\"a}hler TQFT.

\begin{enumerate}
\item The first observable is one we have already mentioned, namely
the inclusion of additional copies of $\Theta$ to bring the degree of
a trivalent graph in $\mathcal{A}(\emptyset)$ up to $k$, for a
hyperk{\"a}hler manifold $X$ of real-dimension
$4k$. More precisely, suppose $D\in\mathcal{A}(\emptyset)$ has degree
$l$, that is $2l$ vertices, where $l<k$. Then
$$W_X(D)\in\mathrm{H}^{2l}(X,\mathcal{O}_X).$$
In order to integrate this, as in the original Rozansky-Witten
invariants, we first need to multiply by
$$[\bar{\o}^{k-l}]\in\mathrm{H}^{2(k-l)}(X,\mathcal{O}_X)$$
or equivalently 
$$W_X(\Theta)^{k-l}\in\mathrm{H}^{2(k-l)}(X,\mathcal{O}_X).$$
More generally, an arbitrary element of
$\mathrm{H}^{2(k-l)}(X,\mathcal{O}_X)$ may be included as an
observable, though for irreducible $X$ this cohomology space is
one-dimensional and hence generated by $W_X(\Theta)^{k-l}$.

\item The second type of observable arises from taking a knot $K$ in
our 3-manifold $M$ and a holomorphic vector bundle $E$ over $X$. We
then take the trace in a fibre $\mathbb E$ of $E$ of the holonomy
around the knot $K$ of the gauge field $A$, to get
$$\mathcal{O}(K;{\mathbb E})=\mathrm{Tr}_{\mathbb E}\exp (\oint_K A).$$
Note that originally Rozansky and Witten introduced a similar
observable with $\mathbb E$ a representation of
$\mathrm{Sp}(k)$. Since $X$ has holonomy $\mathrm{Sp}(k)$ (or
contained in $\mathrm{Sp}(k)$ if $X$ is reducible), the frame bundle
of $X$ is a principal $\mathrm{Sp}(k)$-bundle. Thus the representation
$\mathbb E$ induces a bundle $E$ over $X$, which is a tensor
bundle. Thompson's construction is a generalization to arbitrary
holomorphic\footnote{When $X$ is hyperk{\"a}hler (not just holomorphic
symplectic) we can obtain a `good' observable by adding the condition
that $E$ be hyper-holomorphic (see~\cite{thompson00}). This gives the
observable some invariance properties under variation of the
compatible complex structure on $X$. For example, all tensor bundles
are hyper-holomorphic.} vector bundles $E$ over $X$.

\item The third type of observable comes from choosing a $(0,q)$-form
$\la$ on $X$ with values in
$\La^{l_1}T\otimes\cdots\otimes\La^{l_r}T$, where $b_1(M)=r>0$. See
Thompson~\cite{thompson99} for how these are constructed.
\end{enumerate}

In the present context, the first kind of observable is not relevant
as we are working over the commutative ring $\mathcal H$, and so we do
not require our cohomology classes to be of top degree.

The second kind of observable comes from including a knot $K$ and
associating a holomorphic vector bundle $E$ to it. The inclusion of a
knot in our 3-manifold $M$ means that instead of the LMO invariant of
$M$, the Murakami-Ohtsuki TQFT will give us the Kontsevich integral
$$Z(K\subset M)\in\mathcal{A}(S^1)$$
of the knot $K$ in $M$. If we then apply the hyperka{\"a}hler weight
system to $Z(K\subset M)$, as in Thompson~\cite{thompson00} or the
author's article~\cite{sawon00}, we get the expectation value of the
observable, instead of the usual partition function.

\begin{figure}[htb]
\includegraphics[scale=0.7]{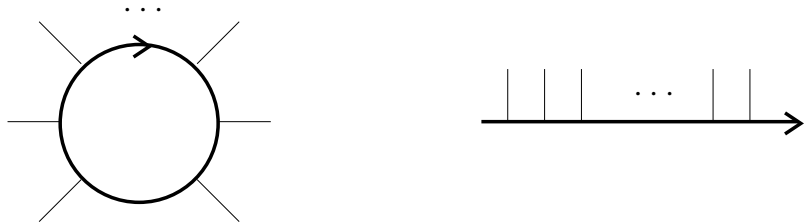}
\caption{A wheel and a comb (or a hedgehog and a flat hedgehog).}
\label{figure9}
\end{figure}

\begin{prop}
The dependence on the vector bundle $E$ is purely through its Chern
character. 
\label{chern}
\end{prop}
\noindent
{\bf Proof:} As with the usual construction of a hyperk{\"a}hler
weight system, when we include a vector bundle $E$ we place its Atiyah
class 
$$\alpha_E\in\mathrm{H}^1(X,T^*\otimes\mathrm{End}E)$$
at the univalent vertices, which lie on the circle $S^1$. We then
contract indices and complete the construction as before. Thus $E$
enters into the construction through the appearance of the cohomology
classes
$$\mathrm{Tr}(\alpha_E^l)\in\mathrm{H}^l(X,\otimes^lT^*)$$
coming from {\em wheels\/} (also known as {\em hedgehogs\/}) with $l$
spokes, as shown in Figure~\ref{figure9}. Note that powers of
$\alpha_E$ are obtained by composing elements of $\mathrm{End}E$ and
taking the cup-product in cohomology. The traces
$\mathrm{Tr}(\alpha_E^l)$ are known as the {\em big Chern classes\/}
of $E$ (see Kapranov~\cite{kapranov99}). 

Recall that we have isomorphisms
$$\mathcal{B}\stackrel{\chi}{\longrightarrow}\mathcal{A}(I)\stackrel{\rho}{\longrightarrow}\mathcal{A}(S^1).$$
By virtue of the isomorphism $\rho$, we can break a wheel to make a
{\em comb\/}, as shown in Figure~\ref{figure9}. By virtue of the
isomorphism $\chi$, we may assume that the legs of an element of
$\mathcal{A}(I)$ (and hence of an element of $\mathcal{A}(S^1)$) are
arranged symmetrically. In the hyperk{\"a}hler context, symmetry
becomes antisymmetry, and in effect we have projected the big Chern
class to the exterior product, where it becomes the usual Chern class 
$$\mathrm{ch}_l(E)\in\mathrm{H}^l(X,\La^lT^*)=\mathrm{H}^{l,l}_{\bar{\partial}}(X)\subset\mathrm{H}^{2l}(X).$$
This completes the proof.\hfill$\Box$

\vspace*{5mm}
These observables can be generalized to links with $l>1$ components,
where we associate a holomorphic vector bundle to each link
component. By generalizing the above argument we can show that the
dependence on these vector bundles is once again purely through their
Chern classes. The crucial point is the symmetry of the diagrams which
ultimately reduces the big Chern class to the usual Chern
class. Recall that we have maps 
$$\mathcal{B}_g\stackrel{\chi}{\longrightarrow}\mathcal{A}(\coprod_{i=1}^gI)\longrightarrow\mathcal{A}(\coprod_{i=1}^gS^1).$$
The first of these is an isomorphism, and implies that the legs of an
element of $\mathcal{A}(\coprod_{i=1}^gI)$ may be arranged in $g$
symmetric collections. The second map, while not an isomorphism, is
surjective. Indeed $\mathcal{A}(\coprod_{i=1}^gS^1)$ is isomorphic to
the quotient of $\mathcal{A}(\coprod_{i=1}^gI)$ by the link relations
(this is Theorem 3 in Bar-Natan et al~\cite{bgrtII99}). This implies
that the legs of an element of $\mathcal{A}(\coprod_{i=1}^gS^1)$ may
also be arranged symmetrically, and the rest of the argument follows
as before.

Let us return to the knot case. In terms of the TQFT we can think of
this observable in the following way. A toroidal neighbourhood $N(K)$
of the knot is a solid torus with an embedded $S^1$. Thought of on its
own, the solid torus is a 3-manifold with boundary a genus one Riemann
surface, and hence gives us an element
$$Z(N(K))\in\mathcal{H}_1=\bigoplus_q\mathrm{H}^q(X,\La^{\bullet}T).$$
The complement of $N(K)$ in $M$ also gives rise to an element
$$Z(M\backslash N(K))\in\mathcal{H}_1=\bigoplus_q\mathrm{H}^q(X,\La^{\bullet}T),$$
and pairing as in the previous section gives
$$Z(M)=\langle Z(N(K)),Z(M\backslash
N(K))\rangle_{\mathcal{H}_1}\in\mathcal{H}=\bigoplus_q\mathrm{H}^q(X,\mathcal{O}_X).$$
Note that the space $\mathcal{H}_1$ is the Hochschild cohomology of
$X$, as on page 67 of Kontsevich~\cite{kontsevich99ii} where it is
denoted 
$$\mathrm{HH}^m(X):=\bigoplus_{p+q=m}\mathrm{H}^q(X,\La^pT).$$
The Hodge cohomology
$$\bigoplus_{i,j}\mathrm{H}^j(X,\La^iT^*)=\bigoplus_{i,j}\mathrm{H}^{i,j}_{\bar{\partial}}(X)\subset\bigoplus_n\mathrm{H}^n(X)$$
of $X$ acts linearly on the Hochschild cohomology $\mathcal{H}_1$,
with the action given by taking the cup product in cohomology combined
with the convolution operator
$$\La^iT^*\otimes\La^pT\longrightarrow\La^{p-i}T$$
acting fiberwise. Pure Hodge classes, i.e.\ those lying in
$$\bigoplus_i\mathrm{H}^i(X,\La^iT^*),$$
preserve the $\mathbb{Z}-$grading on $\mathcal{H}_1$ (though not the
$(p,q)$ bi-grading).

The observable is added by including the embedded $S^1$ with its
associated vector bundle $E$. This changes the element $Z(N(K))$ by
twisting by some element $\d(E)$ of Hodge cohomology, using the
action described above. By virtue of Proposition~\ref{chern}, $\d(E)$
must be a characteristic class of $E$, and presumably it is the Chern
character. In any case, it will be of pure Hodge type and so the twist
preserves the degree of $Z(N(K))$. Pairing the twisted term with
$Z(M\backslash N(K))$ as before gives us the expectation value of the
observable
$$Z(M,\mathcal{O}(K;E))=\langle \d(E).Z(N(K)),Z(M\backslash N(K))\rangle_{\mathcal{H}_1}\in\mathcal{H}.$$
The action of the Hodge cohomology on the Hochschild cohomology is the
hyperk{\"a}hler analogue of a well-known diagrammatic operation and
will be discussed in more detail in a forthcoming paper of Roberts and 
Willerton~\cite{rw00}.

It is not clear how this interpretation can be modified to describe
the third kind of observable, but let us make the following
comments. These observables are defined for 3-manifolds with
$b_1(M)=r>0$, whereas the Murakami-Ohtsuki TQFT is only defined on
rational homology spheres. So we need to assume that the
hyperk{\"a}hler TQFT can be extended to arbitrary
3-manifolds. Now involved in the construction is a choice of basis
$\{\g_i\}$ for $\mathrm{H}_1(M)$ (see Thompson~\cite{thompson99})
which could be regarded as a $b_1(M)=r$ component link. If we consider
a neighbourhood of the link, as we did for knots earlier, we are led
to 3-manifolds with boundaries consisting of $r$ connected components,
each component a genus one Riemann surface. Of course we don't yet
have a definition of the TQFT for such 3-manifolds when
$r>1$. Nonetheless, the form $\la$ should be holomorphic and hence we
have a Dolbeault cohomology class 
$$[\la]\in\bigoplus_q\mathrm{H}^{0,q}_{\bar{\partial}}(X,\La^{l_1}T\otimes\cdots\otimes\La^{l_r}T)=\bigoplus_q\mathrm{H}^q(X,\La^{l_1}T\otimes\cdots\otimes\La^{l_r}T)\subset\mathcal{H}_r$$
which sits inside the Hilbert space associated to a genus $r$ Riemann
surface. This suggests that $\mathcal{H}_r$ must be closely related to
the Hilbert space associated to $r$ disjoint genus one Riemann
surfaces: the latter is probably a quotient space of the former.

We can be a little more specific when $r=1$. Then the generator $\g_1$
of $\mathrm{H}_1(M)$ can be regarded as a knot $K$ in $M$, and as
before we get 
$$Z(N(K))\qquad\mbox{and}\qquad Z(M\backslash N(K))\in\mathcal{H}_1.$$
We can discard $Z(N(K))$ and pair $Z(M\backslash N(K))$ with
$[\la]\in\mathcal{H}_1$ instead, giving 
$$\langle [\la],Z(M\backslash
N(K))\rangle_{\mathcal{H}_1}\in\mathcal{H}$$
which presumably is the expectation value of the observable.

Note that in general given a 3-dimensional TQFT we can construct knot
invariants in the following way. By removing a toroidal neighbourhood
of the knot we get a 3-manifold with boundary a torus. Applying the
TQFT we therefore get a vector in the Hilbert space associated to the
torus. To get knot invariants we can take the components of this
vector with respect to some choice of basis of the Hilbert
space. Equivalently, we can pair this vector with another vector in
the Hilbert space. In effect all the observables in this section can
be interpreted in this way.

\section{Appendix}

We now give the proof of Proposition~\ref{rho} which we postponed
earlier on. In fact the result we shall prove is slightly more
general, and is due to Lescop~\cite{lescop00}. The author is grateful
to Christine Lescop for devising this proof and explaining it to him.

\begin{figure}[htb]
\includegraphics[scale=0.7]{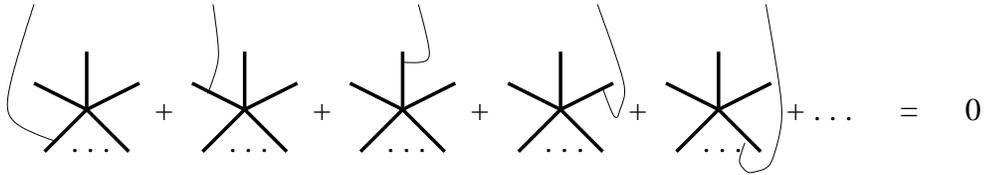}
\caption{Branching relations for higher valency vertices.}
\label{figureA1}
\end{figure}

Proposition~\ref{rho} says that the space of chord diagrams on the
collection of $g$ intervals is isomorphic to the space of chord
diagrams on the chain graph $\G_g$, that is
$$\rho
:\mathcal{A}(\coprod_{i=1}^gI)\stackrel{\cong}{\longrightarrow}\mathcal{A}(\G_g).$$
The map $\rho$ is given by attaching a certain (unitrivalent) tree to
the intervals. More generally, let $\D$ be an arbitrary tree which has
precisely $2g$ univalent vertices. By `arbitrary' we mean that some of
the vertices may have valency greater than three. When we attach this
to the collection of $g$ intervals (in some way), the resulting graph
$\G$ may also have vertices of valency greater than three. We can
consider the space of chord diagrams on $\G$ modulo the AS, IHX, STU,
and branching relations, though we need to add additional branching
relations for vertices of higher valency (see
Figure~\ref{figureA1}). We denote the resulting space
$\mathcal{A}(\G)$ as before.

\begin{figure}[htb]
\includegraphics[scale=0.7]{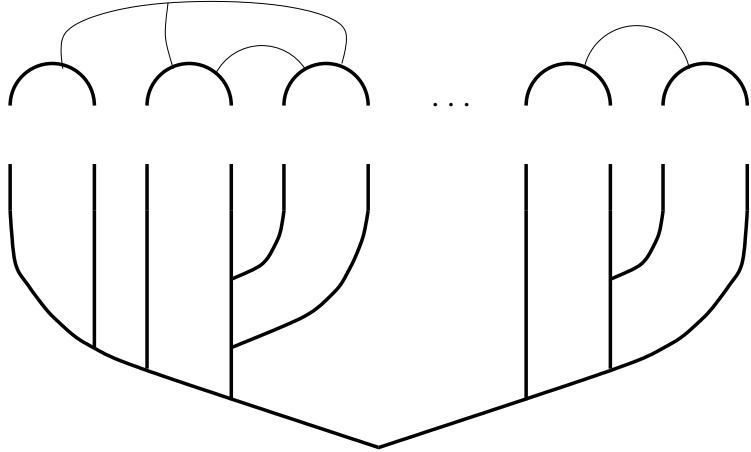}
\caption{Joining a tree to a collection of intervals.}
\label{figureA2}
\end{figure}

The operation of attaching the tree $\D$ to the collection of $g$
intervals as shown in Figure~\ref{figureA2} results in a map
$$\rho :\mathcal{A}(\coprod_{i=1}^gI)\longrightarrow\mathcal{A}(\G).$$
We use the same notation as in Proposition~\ref{rho} as this is a
generalization of that map.

By constructing an explicit inverse to $\rho$ we shall prove the
following.

\begin{figure}[htb]
\includegraphics[scale=0.7]{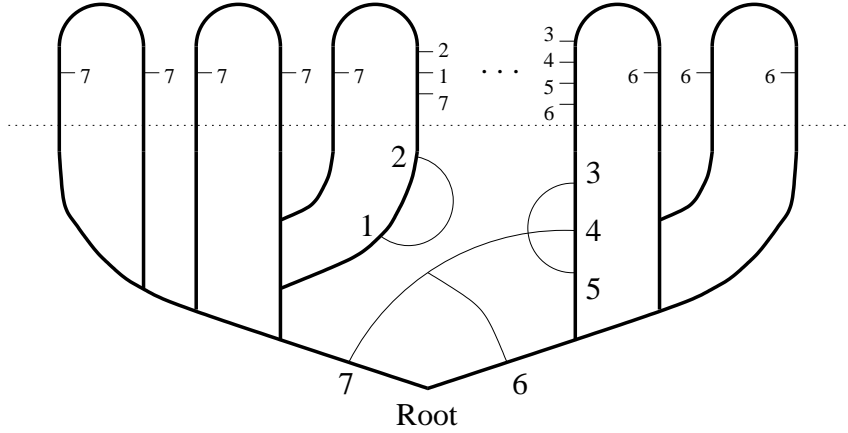}
\caption{Pushing a chord diagram off a tree and onto a collection of intervals.}
\label{figureA3}
\end{figure}

\begin{prop}
The map $\rho$ is an isomorphism of $\mathcal{A}(\emptyset)$-modules.
\end{prop}
\noindent
{\bf Proof:} Recall that a chord diagram $D$ on $\G$ is the union of
$\G$ and a unitrivalent graph, the chord graph $Q$. Given $D$, we want
to `push' all the legs of the chord graph $Q$ off the `tree part' of
$\G$ and onto the `interval part' of $\G$. In Figure~\ref{figureA3}
this corresponds to pushing the legs up past the dotted line. The
branching relations allow us to push legs past the vertices of
$\G$. In order to turn this into a rigorous proof, we need to make
this `pushing' operation more precise. The following idea is due to
Lescop~\cite{lescop00}.

We may assume that we begin with no legs above the dotted line. We
label all the legs of $Q$ by natural numbers $1$, $2$, $3$, etc. Next
we choose a generic point of $\G$ which lies below the dotted line,
and call it the {\em root\/}. The dotted line intersects $\G$ in $2g$
points. If $p$ is such a point, we trace down along $\G$ to the root
making a note of which legs are encountered along the way; then we add
`labelled tabs' to $\G$ above $p$ corresponding to the legs
encountered and {\em in the same order\/}. We do this for all such
points $p$. Figure~\ref{figureA3} shows an example of how the
resulting picture may look.

We are now in a position to define the inverse to $\rho$. Given a
chord diagram $D$ on $\G$, we detach $Q$ from $\G$ but keep note of
the labellings of the legs. Then we take the sum over all the ways to
reattach $Q$ to the labelled tabs above the dotted line. Of course a
leg labelled $i$ may only be attached to a tab labelled $i$. Since
there are no longer any legs attached to the part of $\G$ below the
dotted line (that is, the tree part $\D$), we may remove this part and
are left with a sum of chord diagrams on a collection of $g$
intervals. We call this $\s(D)$.

As an element of $\mathcal{A}(\G)$, we only regard $D$ up to the AS,
IHX, and STU relations. It is easy to see that if $D$ and $D^{\prime}$
are equivalent under these relations, then $\s(D)$ and
$\s(D^{\prime})$ are also equivalent under the AS, IHX, and STU
relations (applied to chord diagrams on the collection of $g$
intervals). It is also clear from the way $\s$ is defined that if $D$
and $D^{\prime}$ are equivalent under the branching relations ({\em
not\/} including the branching relation at the root) then
$\s(D)=\s(D^{\prime})$. This means that the only possible ambiguity in
the construction concerns the root.

\begin{figure}[htb]
\includegraphics[scale=0.65]{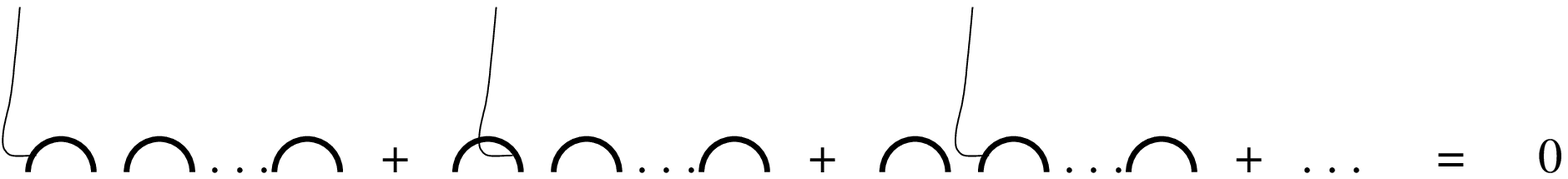}
\caption{Relation in $\mathcal{A}(\coprod_{i=1}^gI)$.}
\label{figureA4}
\end{figure}

If we slide a leg of $D$ past the root to get $D^{\prime}$ then the
difference between $\s(D)$ and $\s(D^{\prime})$ will be a sum of terms
as on the left hand side of the equation in
Figure~\ref{figureA4}. Note that we have only drawn the legs which are
closest to the endpoints of the intervals, and there may be many more
legs which we have not shown.

However, the vanishing of this sum is a well-known relation in the
space $\mathcal{A}(\coprod_{i=1}^gI)$. Thus $\s(D)$ and 
$\s(D^{\prime})$ are equivalent chord diagrams. If the root happens to
be a vertex of $\G$, then `sliding a leg past the root' would mean
applying the branching relation at that vertex, and the conclusion is
still valid. Instead of moving legs past the root, we could move the
root itself, thus resulting in a new map $\s^{\prime}$. An argument
analogous to the one above shows that $\s(D)$ and $\s^{\prime}(D)$ are
equivalent chord diagrams.

Therefore we have a well-defined map
$$\s :\mathcal{A}(\G)\longrightarrow\mathcal{A}(\coprod_{i=1}^gI).$$
It is now a simple exercise to verify that $\s$ and $\rho$ are
inverses, thus proving the proposition. \hfill $\Box$

\end{document}